\documentclass{jdcs}
\usepackage[active]{srcltx}
\usepackage[T1]{fontenc}
\usepackage[utf8]{inputenc}
\usepackage{lmodern}
\usepackage[a4paper]{geometry}
\usepackage{amsmath,amssymb,amsthm,mathabx}
\usepackage{pgf,tikz}
\usepackage{mathrsfs}
\usepackage{graphicx}
\usepackage{xcolor}
\usepackage[colorlinks=true,breaklinks=true,linkcolor=blue,urlcolor=green,citecolor=red]{hyperref} 
\usetikzlibrary{arrows}

\newcommand{\R}{\mathbb{R}}

\renewcommand{\epsilon}{\varepsilon}
\renewcommand{\d}{\mathrm{d}}
\newcommand{\p}{\partial}
\newtheorem{theorem}{Theorem}
\newtheorem{lemma}{Lemma}
\newtheorem{proposition}{Proposition}

\newtheorem{remark}{Remark}
\newtheorem{definition}{Definition}
\theoremstyle{remark}

\def\Fin#1{\leavevmode\unskip\nobreak\quad\hspace*{\fill}{#1}}
\begin{document}
\title{Singularities of optimal time affine control systems: the limit case}
\author{Orieux, M.}\address{Scuola Internazionale Superiore di Studi Avanzati
 "SISSA", Via Bonomea, 265, 34136 Trieste TS Italy}
\author{Roussarie, R.}\address{Institut de Mathématiques de Bourgogne, UMR 5584, Universi\'e de Bourgogne, 9 avenue Alain Savary
BP 47870
21078 Dijon, France}

\date{15/10/2020}
\begin{abstract} We study the singularities of the extremal flow - and thus, of the exponential map - for minimum time control-affine problems in 4D with 2D controls. The case of mechanical systems has been previously studied in \cite{orieux:2018:singularities}. Here we provide a unifying view point with a parameter and close the final case remaining open, proving the flow is piece-wise smooth. After regularization, the problem boils down to the study of a bifurcation around some nilpotent equilibrium in the singular locus.
\end{abstract}

\maketitle
\section*{Introduction}
Control-affine systems generalize sub-Riemannian geometry by adding a drift, 
and arise naturally from controlled mechanical systems. We handle the case 
of minimizing the final
time, for an affine control system defined on a connected $4$-dimensional 
manifold $M$ under a generic assumption given
below. The control $u$ is in dimension $2$ (double input case), and contained in the
unitary euclidean ball, $B$:
$$\dot{x}=F_0(x)+u_1F_1(x)+u_2F_2(x),\; x\in M, u\in B$$
where the $F_i$'s are smooth vector fields on $M$.
In this case, the research of optimal
trajectories leads, according to Pontrjagin's Maximum Principle, to study a singular Hamiltonian system
defined on the cotangent bundle of the manifold $M$. We are interested in the 
local behavior of the flow of this Hamiltonian system, called extremal flow, around its
singularities. The singular locus, also called
the switching set, is a codimension two submanifold where a switch - 
a discontinuity of the optimal control - is susceptible
to occur. We compare this case with single-input control systems where singularities are of codimension one, i.e.,
when the control is scalar:
$\dot{x}=F_0(x)+uF_1(x)$. 
They have been extensively studied, and many things are known regarding their singularities, 
see for instance the monograph \cite{intro}.
In the present paper, we give a precise description of the behavior of
the extremal flow around the singular locus. 
We answer two important questions which remained open: 

- Is there trajectories crossing each point of the singular locus ?

- What is the regularity of the lifted flow in that case ? 

\noindent It is known since \cite{caillau2012minimum} that the singular locus can be partitioned into three subsets, 
giving three different configurations for the flow. The first two cases where handle in \cite{biolo} and \cite{orieux:2018:singularities}. The last one, detailed below, was more difficult to handle,
being the frontier of a bifurcation phenomenon, and a higher order analysis of the dynamics was necessary.
Those questions have simpler answers when one consider the problem with the control set $U$ 
being a polyhedron (which makes no difference in the single input case), and have been treated for the minimum time case in \cite{controlbox}, for instance.
 We attempt to give a unify point of view for those three cases,
through a formulation with a parameter.

	\section{Setting}
	Let $M$ be a 4-dimensional manifold, $x_0$, $x_f\in M$ and consider the following time
	optimal control system:
	\begin{equation}
	\begin{cases}
	\dot{x}(t)=F_0(x(t))+u_1(t)F_1(x(t))+u_2(t)F_2(x(t)), \quad
	t\in[0, t_f], \quad u \in U\\
	x(0)=x_0 \\
	x(t_f)=x_f\\
	t_f\rightarrow \min.
	\end{cases}
	\label{Tmin}
	\end{equation}
	 where the $F_i$'s are smooth vector
	fields. Set $F_{ij}:=[F_i,F_j]$, and $H_{ij}=\{H_i,H_j\}$. The pont of view of this paper is local, let's pick a point $\bar{x}\in M$. We make the following generic assumption on the distribution $\mathcal{D}=\{F_0,F_1,F_2\}$.
	\begin{equation}
	\det(F_1(\bar{x}), F_2(\bar{x}), F_{01}(\bar{x}), F_{02}(\bar{x}))\neq0.
	\label{eq:A}
	\tag{\textcolor{blue}{A}}
	\end{equation}
	This is slightly stronger than the assumption that the distribution
	$\mathcal{D}$ is of step 2, but it is natural in the applications 
	to mechanical systems, for instance. 
	By Pontrjagin’s Maximum Principle, optimal trajectories are the projection of integral curves of the
	maximized Hamiltonian system 
	defined on the cotangent bundle of $M$ by
	\[H^{\max}(z):=H_0(z)+\sqrt{H_1^2(z)+H_2^2(z)},\ z\in T^*M\]
	where we have denoted $H_i(z):=\langle p,F_i(x)\rangle$ in canonical coordinate $(x,p)\in T^*_xM$.
	See \cite{orieux:2018:singularities} or \cite{agrachev04} for more details about 
	Pontrjagin's Maximum Principle.
	Besides, it implies the feedback control \[u=\frac{1}{\sqrt{H_1^2+H_2^2}}(H_1,H_2)\]
	whenever $(H_1,H_2)\neq (0,0)$.
	An integral curve of $H^{\max}$ is called an extremal. 
	The extremal system is smooth outside the singular locus, or switching surface, 
	defined by $$\Sigma:=\{H_1=H_2=0\}.$$
	
	\begin{definition}[Bang and bang-bang and singular extremals]
    An extremal $z(t)$ is said to be \emph{bang} if $(H_1,H_2)(z(t))\neq (0,0)$ for all $t$. 
    It is \emph{bang-bang} if it is a concatenation of bang arcs. We say it is \emph{singular} if it is contained in the singular locus.
	\end{definition}
	
	For single input systems, if one consider bounded controls, taking values in $[-1,1]$, then the maximized Hamiltonian is
	$$H(x,p)=H_0(x,p)+|H_1(x,p)|$$ with as before $H_i(x,p)=\langle p,F_i(x)\rangle$, $i=0,1$, and $u=$sign$(H_1)$ when $H_1$ is non-zero. The singular locus for this problem is $\{H_1=0\}$. The singular controls can then be easily calculated by differentiating the relation
	$H_1(z(t))=0$. 
	For the so-called order one singular extremals, we have 
	\begin{proposition}
		The singular flow for single input systems is given by the Hamiltonian $H^s=H_0-\frac{H_{100}}{H_{101}}H_1$, with the singular control given by
		$u^s=\frac{\{\{H_1,H_0\},H_0\}}{\{\{H_1,H_0\},H_1\}}$.
	\end{proposition}
	Discontinuities of the control $u$ along an extremal are called switchings, a time $\bar{t}$ 
	at which a switch occurs is called switching time, and $z(\bar{t})$, a switching point. Bang extremals are 
	the one that do not cross $\Sigma$. 
	
	Now we tackle the double-input system.
	A singular minimum time extremal is such that $H_{12}(z(t))\neq0$, 
	see remark \ref{B} below, one can also see \cite{caillau2012minimum}.
	The following proposition has been proven in \cite{caillau2012minimum}.
	
\begin{proposition}
	There exists a singular flow inside $\Sigma$, on
	which we have the control feedback:
	$u_s=\frac{1}{H_{12}}(-H_{02},H_{01})$, and the singular flow is
	smooth. It is solution of the Hamiltonian system given by
	$\tilde{H}=H_0-\frac{H_{02}}{H_{12}}H_1+\frac{H_{01}}{H_{12}}H_2$.
	\label{sing}
\end{proposition}
	\noindent\textbf{Proof.} The proposition is obtained by differentiating the identically zero switching function
	$(H_1,H_2)(z(t))$ with respect to the time. \Fin{$\Box$}
	
	From \cite{caillau2012minimum} and \cite{orieux:2018:singularities}, we know 
	$\Sigma$ is partitioned into three subsets, leading to three very different local dynamics 
	in their neighborhoods, namely (we use the notation $H_{ij}=\{H_i,H_j\}$)
	\begin{align*}
	\Sigma_- &=\{H_{12}(z)^2<H_{02}(z)^2+H_{01}(z)^2\}\\
	\Sigma_+ &=\{H_{12}(z)^2>H_{02}(z)^2+H_{01}(z)^2\}\\
	\Sigma_0 &=\{H_{12}(z)^2=H_{02}(z)^2+H_{01}(z)^2\}.
	\end{align*} 
	The behavior of the flow in a neighborhood of $\Sigma_0$ remains open,
	 as the case $\Sigma_-$ and $\Sigma_+$ were settled in \cite{orieux:2018:singularities,biolo}, 
	 but we attempt to provide in the next section a unification of the different viewpoints.

\begin{remark}
Note that, since along a Pontrjagin extremal, the adjoint state $p$ cannot vanish, $(A)$ is equivalent to 
$$(H_1^2+H_2^2+H_{01}^2+H_{02}^2)(\bar{z})>0,$$
at any point $\bar{z}=(\bar{x},p)$, $p\in T_{\bar{x}}M$.
\label{B}
\end{remark}
Let us take a look to the single input case. The study of switchings is very simplified by the following
fact: under generic hypothesis, one can define switching times by the implicit function theorem for
order one switchs (switching occur when $H_1$ vanishes, and the Hamiltonian 
vector field $X_{H_0}+X_{H_1}$  is well defined). This fact remain true when $U$ is a box 
instead of a ball, see \cite{controlbox}. The components of the control $u$ just go from $+1$ to $-1$ or the opposite.

\section{Formulation with a parameter} 
In this section we introduce rather artificially a parameter in the previous
 dynamical system in order to unify the viewpoints.
 Thanks to (\ref{eq:A}), one can make the
change of coordinates: 
\[z=(x,p)\in T^*M\mapsto (x,H_1,H_2,H_{01},H_{02}) \in M\times\mathbb{R}^4.\]
Then one uses a polar blow up by setting $(H_1,H_2)=\rho(\cos\theta,\sin\theta)$ 
and \linebreak
 $(H_{01},H_{02})=r(\cos\phi,\sin\phi)$. 
The dynamics boils down to the system:
\begin{equation} 
\begin{cases}
\dot{\rho}=r\cos(\theta-\phi)\\
\dot{\theta}=\frac{1}{\rho}(H_{12}-r\sin(\theta-\phi))\\
\dot{\xi}= h(\rho,\theta,\xi)
\end{cases}
\label{syst1}
\end{equation} 
where $\xi=(x,r,\phi)$ and $h$ is a smooth function defined on an open set $O$ of 
$\mathbb{R}\times\mathbb{R}\times D$, $D$ being a compact domain of $\mathbb{R}^6$; $h$ 
has values in $\mathbb{R}^6$. We set $\psi=\theta-\phi$, and rescale the time according to 
$dt_1=rdt$: by remark \ref{B}, $r$ is never $0$ in a neighborhood of $\Sigma$, 
meaning this defines a diffeomorphism and new dynamics is conjugate the one of 
system \ref{syst1}. This boils down to study a general system with the following structure 
(the derivation with respect to the time $t_1$ still being noted "$^.$"):  
\begin{equation}\begin{cases}
\dot{\rho}=\cos\psi \\
\dot{\psi}=\frac{1}{\rho}(g(\rho,\psi,\xi)-\sin\psi)\\
\dot{\xi}=\tilde{h}(\rho,\psi,\xi)
\end{cases}
\label{d}
\end{equation}
where $g$, $\tilde{h}$ are smooth functions
defined on an open set $O$ of $\mathbb{R}\times\mathbb{R}\times D$, $D$ being a
compact domain of $\mathbb{R}^k$, for any $k\geq1$ (in our problem $k=6$); $g$ depends smoothly on $\rho(\cos\theta,\sin\theta)$. 
One can set $s=\psi-\pi/2$. By a small abuse of the notations, we still note 
$g(\rho,s,\xi)=g(\rho,s+\pi/2,\xi)$, and we have
$g(\rho,s,\xi)=a(\xi)+O(\rho)$ near $\rho=0$. 

 The three cases $\Sigma_+$, $\Sigma_-$ and $\Sigma_0$ 
correspond to the values $g(\bar{z})>1$, $g(\bar{z})<1$ 
and $g(\bar{z})=1$, with $\bar{z}=(0,\bar{s},\bar{\xi})\in \Sigma$, see \cite{orieux:2018:singularities} for more details. 
We set $a(\xi)=1+\alpha+a_0(\xi)$, with $a_0(\bar{\xi})=0$. So that the bifurcation parameter 
giving the three cases is $\alpha=a(\bar{\xi})-1$ and 
$$g(\rho,s,\xi)=1+\alpha+a_0(\xi)+O(\rho).$$
Then (\ref{d}) in the new time becomes (with slight abuse of notation):
\begin{equation}
(Z): \begin{cases}
\dot{\rho}=-\sin s \\
\dot{s}=\frac{1}{\rho}(1+\alpha+a_0(\xi)-\cos s+O(\rho))=\frac{G_\alpha(\rho,s,\xi)}{\rho}\\
\dot{\xi}=\tilde{h}(\rho,s,\xi)
\end{cases}
\label{e}
\end{equation}

\begin{remark}
We actually have $a(x,H_{01},H_{02})=\frac{H_{12}}{\sqrt{H_{01}^2+H_{02}^2}}(x,0,0,H_{01},H_{02})$. 
\end{remark} 
\subsection{The case $\Sigma_-$.}

In the system above, $\bar{z}\in\Sigma_-$ if and only if $\alpha<0$. 
That case was settled in \cite{orieux:2018:singularities} by theorem 
\ref{theorem} recalled below:

\begin{theorem}
 In a neighborhood $O_{\bar{z}}$ with $\bar{z}\in\Sigma_-$, existence and uniqueness hold, all extremal are bang-bang, with at most 
 one switch. The extremal flow $z:(t,z_0)\in[0,t_f]\times O_{\bar{z}}\mapsto z(t,z_0)\in M$ is
	piecewise smooth. More precisely, $O_{\bar{z}}$ can be stratified as follows:
	$$O_{\bar{z}}=S_0\cup S^s\cup\Sigma$$ 
	where $S^s$ is the codimension-one submanifold of initial
	conditions leading to the switching surface,
	$S_0=O_{\bar{z}}\setminus (S^s\cup \Sigma)$.
	Both are stable by the flow, which is smooth on $[0,t_f]\times S_0$, and 
	on $[0,t_f]\times S^s\setminus\Delta$  where $\Delta=\{(t_\Sigma(z_0),z_0),\; z_0\in S^s\}$, 
	and $t_\Sigma(z_0)$ is the switching time of the extremal initializing at $z_0$, 
	and continuous on $O_{\bar{z}}$.
	\label{theorem}
\end{theorem}

In \cite{orieux:2018:singularities}, the authors also studied the regular-singular 
transition between the strata, and exhibited $\log$-type singularities.
\begin{remark}
By proposition \ref{sing}, there is no admissible singular flow contained in $\Sigma_-$,
otherwise $\|u_s\|^2=\frac{H_{02}^2+H_{01}^2}{H_{12}^2}>1$ which violates our
constraint on $u$. 
\end{remark}

\subsection{The case $\Sigma_+$.}

This corresponds to $\alpha>0$. We prove the following, see also \cite{biolo}, theorem 3.5: 
\begin{proposition}
	In a neighborhood of a point $\bar{z}$ in $\Sigma_+$, there is no switch, and the extremal flow is smooth,
	 i.e., $\Sigma_+$ is never crossed. In other words, $\rho$ does not vanish in $(\ref{e})$.
\end{proposition}

\noindent\textbf{Proof.} By the analysis above, this boils down to prove 
that, if $\alpha>0$, along an extremal $z$, $\rho$ never vanishes in a (relatively compact) neighborhood $\bar{O}$ of $(0,0,0)\in\mathbb{R}\times\mathbb{R}\times\mathbb{R}^k$. 
Set $f$ such that  $\rho\dot{s}=1+\alpha+f(\rho,s,\xi)-\cos s:=\Theta$  in such a neighborhood, the differential of $f$ is bounded by below by a negative constant $-a<0$. 
$$\frac{d}{dt}(\rho\Theta)=\dot{\rho}(1+\alpha-\cos s+f(z))+\rho(\sin s\dot{s}+df(z).\dot{z})=\rho df(z).Z(z)$$
We get $\frac{d}{dt}(\rho\Theta)>-a\rho=\rho\Theta(-a/\Theta)$. Eventually, since $\alpha>0$ and $u(0,0,0)=0$, 
if $\bar{O}$ is small enough, there exists two positive constant $K$, $k$ with $K>\Theta>k>0$. 
So that, along an extremal
\[\frac{d}{dt}(\rho\Theta)>-\frac{a}{k}\rho\Theta.\] 
By integration on a arbitrary time interval $[0,t]$, we end up with 
\[\rho(t)>\frac{\rho_0\Theta_0}{K}e^{-\frac{a}{k}t},\]
and the proposition follows. \Fin{$\Box$}
\\

Despite the absence of switch, there exists a singular flow inside $\Sigma_+$, however, singular extremals lying in $\Sigma_+$ cannot be 
optimal by the Goh condition, \cite{caillau2012minimum}.
 
\subsection{The bifurcation $\alpha=0$: case $\Sigma_0$.}

This is the main topic of this paper. 
In this case, the two equilibria considered in $\Sigma_-$ merge, 
and we obtain one nilpotent equilibrium that needs desingularization.
Nevertheless - set  $f(x,H_1,H_2,H_{01},H_{02})=\frac{H_{12}}{\sqrt{H_1^2+H_2^2}}$ - under the generic condition

\begin{equation}
\begin{split}
\frac{\p f}{\p x}(\bar{z}) .
(F_0(\bar{x})-\sin\bar{\phi}F_1(\bar{x})+\cos\bar{\phi}F_2(\bar{x}))+\frac{\p f}{\p H_{01}}(\bar{z}) (H_{001}(\bar{z})+H_{101}(\bar{z}))\\
+\frac{\p f}{\p H_{02}}(\bar{z}) (H_{002}(\bar{z})+H_{102}(\bar{z}))\neq0
\end{split}
\label{generica}
\end{equation}
where $\bar{\phi}=\arg(H_{01},H_{02})(\bar{z}),$
we will prove
\begin{theorem}
	For generic systems (\ref{Tmin}), meaning, if assumption (\ref{eq:A}) and (\ref{generica}) holds:
	Let $\bar{z}$ be in $\Sigma_0$ either: there exists a unique trajectory passing 
	through $\bar{z}$, or there exist unique trajectory going out of $\Sigma_0$ at $\bar{z}$. 
	\label{switchnilp}
\end{theorem}
\noindent This result contradicts the last part of theorem 3.5, in \cite{biolo}, 
a counter example in a particular case was given in \cite{dario} (nilpotent case), another one is provided at the end of this paper. Figure 1
is a sketch of the phase portrait.

\begin{remark}
   In the first case, the extremal can be connected to the singular flow in $\Sigma_0$.  
\end{remark}
\begin{figure}
	\centering
	\includegraphics[width=0.8\linewidth]{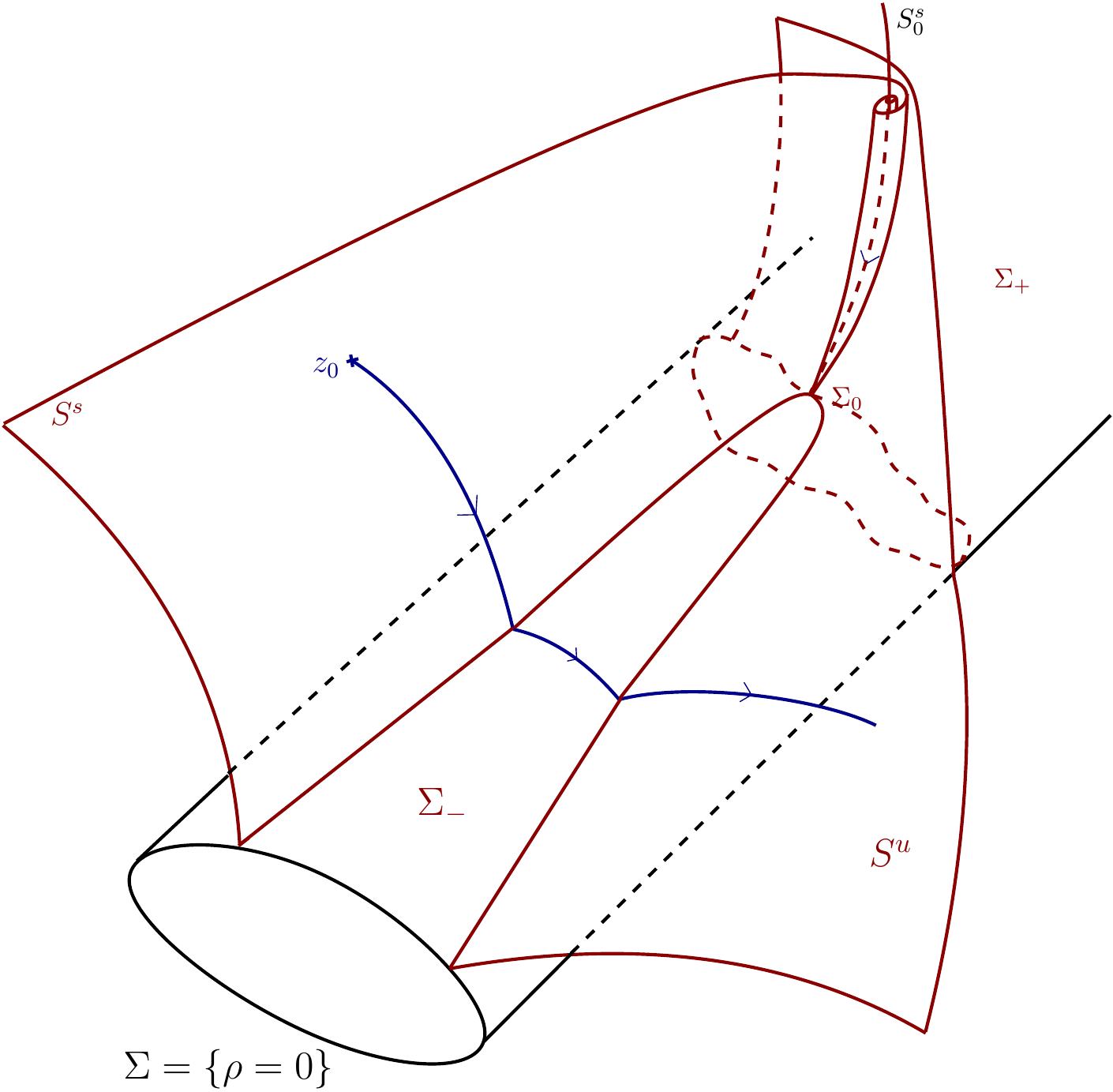}
	\caption{The stable and unstable manifold of $\Sigma_-$ merging on $\Sigma_0$.}
\label{fig0}
\end{figure}
\noindent The regularity of the extremal flow is of theoretical and numerical importance, and we have:
\begin{theorem}
	In a neighborhood $O_{\bar{z}}$ of a point $\bar{z}\in\Sigma_0$, the flow is well defined, and continuous. 
	
	More precisely, there exists a stratification of $O_{\bar{z}}$ into smooth submanifolds:
	$$O_{\bar{z}}=S_0\cup S^s\cup S^s_0$$ 
	where
	\begin{itemize}
\item[--] $S^s_0$ is the submanifold of codimension 2 of initial conditions leading to $\Sigma_0$,
\item[--] $S^s$ is the submanifold of codimension 1 of initial conditions leading to $\Sigma_-$,
\item[--] $S_0=O_{\bar{z}}\setminus(S^s_0\cup S^s)$.
	\end{itemize}
	The extremal flow is smooth restricted to $S_0\times[0,t_f]$, $S^s\times[0,t_f]\setminus\Delta$ and $S^s\times[0,t_f]\setminus\Delta_0$, with $\Delta_0=\{t_{\Sigma_0}(z_0),z_0),\;
	z_0\in S^s_{0}\}$.
	\label{flotnilpreg}
\end{theorem}
Namely, the previous theorem states that the extremal flow is smooth restricted
to each strata, except, obviously, at the point and time where the singular locus
is crossed (at those points, it is only Lipschitz in time alone). 

In the process we also obtain the jumps occurring on the extremal control, a
switching on the control is called a $\pi$-singularity if it is a instant rotation of angle $\pi$, see \cite{orieux:2018:singularities}.
\begin{proposition}
	Consider the extremal $z(t)$ entering the singular locus in $z(\bar{t})=\bar{z}\in\Sigma_0$,
	\begin{itemize}
		\item[-] If $H_{12}(\bar{z})=r(\bar{z})$, the extremal control is continuous,
		\item[-] If $H_{12}(\bar{z})=-r(\bar{z})$, the extremal control has a $\pi$-singularity at time $\bar{t}$.
	\end{itemize} 
	\label{controlcont}
\end{proposition}

\subsubsection{Proof of theorem \ref{switchnilp}}
To that end, we give a precise 
description of the phase portrait of the flow around such an equilibrium. 
 Make the following 
change of 
time $dt_1=\rho dt_2$ to regularize the vector field $Z$, and denote $'=\frac{\d}{\d t_2}$ to get ($\alpha=0$):
\begin{equation}
\begin{cases}
\rho'=-\rho\sin s \\
s'=1+a_0(\xi)-\cos s+O(\rho)=G_0(\rho,s,\xi)\\
\xi'=\rho\tilde{h}(\rho,s,\xi)
\end{cases}
\label{0}
\end{equation}
Since $\xi$ is constant when $\rho=0$, for the sake of clarity, we will work in a neighborhood of $\bar{\xi}=0$,
 keeping in mind that we don't lose generality (and the results holds for any point).
Assumption (\ref{generica}) is equivalent to 
\begin{equation}
da(0).\tilde{h}(0,0,0)\neq 0.
\label{generic}
\end{equation}
Then, we can order the coordinates of 
$\xi=(\xi_1,\xi_2\dots,\xi_k)$ such that $\frac{\p a}{\p \xi_1}(0)\neq 0$. 
This implies $\Gamma=\{G_0(0,s,\xi)=0\}$ is a dimension $k$ manifold 
around $(s,\xi)=(0,0)\in S^1\times\mathbb{R}^k$. We can then chose 
coordinates $\tilde{\xi}=(\tilde{\xi}_1,\dots,\tilde{\xi}_k)=\Phi(\xi)$ 
such that $\tilde{\xi}_1=a_0(\xi)$ meaning, $\Gamma=\{\tilde{\xi}_1+1-\cos s=0\}$. 
Then, set $\zeta:=\tilde{\xi}_1$ to simplify the notations. We obtain

\begin{equation}
\begin{cases}
\rho'=-\rho s+O(\rho s^3)\\
s'=\zeta+s^2/2+O(\rho)+O(|s|^4)\\
\zeta'=c\rho+\rho O(\rho+|s|+|\tilde{\xi}|).
\end{cases}
\label{hot}
\end{equation}
We do not write the dynamics of the other components of $\tilde{\xi}$. They do not
 influence the dynamics of $(\rho,s)$ as we explain below in the paragraph "Back to the original system". Actually, as we will exhibit below, only the first order terms in the derivative of $\zeta$ are relevant for the local dynamics around $0$.
 In the equation $(\ref{hot})$, $c=\tilde{h}(0,0,0)_1$, so that assumption (\ref{generic}) prevent it from
 being $0$.
It will be clear from what follows that the terms of higher order are useless for the local analysis. On $\{\rho=0\}$, the field has two lines of zero with a quadratic contact with $\{\zeta=0\}$. Furthermore the circles $\{\rho=cst\}$ are tangent to the vector field on $\rho=0$. Those lines are normally hyperbolic except at $(\rho,s)=(0,0)$, which is a nilpotent equilibrium.

\paragraph{Blow up.} To study the nilpotent equilibrium $(\rho,s,\zeta)=(0,0,0)$, we will use a specific blow-up
process, called quasi homogeneous blow-up, see \cite{Durmotier} chap. 1: 
$$\begin{cases}
\rho=R^3\bar{\rho}\\ 
s=R\bar{s}\\ 
\zeta=R^2\bar{\zeta}
\end{cases}$$ with $(\bar{\rho},\bar{s},\bar{\zeta})\in S^2_+$ the hemisphere $\rho\geq0$, $R\in\mathbb{R}_+$.
 We will study the dynamics in the two following charts given a vector field smoothly equivalent to the global one we could obtain on $\mathbb{R}\times S^2_+$): 
\begin{itemize} 
\item[$(i)$] For the interior of $S^2_+$, $\bar{\rho}=1$, $(\bar{s},\bar{\zeta})$ in a disc $D^2$, $R\geq0$ in a
 neighborhood of the critical locus of the blow up $R=0$. 
\item[$(ii)$] For the the boundary of $S^2_+$, $(\bar{s},\bar{\zeta})\in \mathbb{S}^1$, in a neighborhood of $(\bar{\rho},R)=(0,0)$. 
\end{itemize}

\subparagraph{The charts $(i)$.} Let us write the dynamics in the blown up coordinates $\varphi(\rho,s,\zeta)=(R^3,R\bar{s},R^2\bar{\zeta})$.
 The blown up vector field $\bar{X}=\frac{1}{R}\varphi_*X$ writes
\begin{equation}
\bar{X}:\begin{cases}
R'=-\frac{1}{3}R\bar{s}+O(R^2)\\
\bar{s}'=\frac{5}{6}\bar{s}^2+\bar{\zeta}+O(R)\\
\bar{\zeta}'=\frac{2}{3}\bar{s}\bar{\zeta}+c+O(R).
\end{cases}
\label{charti}
\end{equation}
Thus, there is a unique equilibrium depending on $c$, which is solution of 
\begin{equation}
\begin{cases}
R=0\\
\bar{\zeta}+\frac{5}{6}\bar{s}^2=0\\
\frac{2}{3}\bar{s}\bar{\zeta}+c=0
\end{cases}
\end{equation}
i.e., $m_0=(0,\bar{s}_0,\bar{\zeta}_0)=(0,\text{sign}(c)(\frac{9}{5}|c|)^{1/3},-\text{sign}(c)\frac{5}{6}(\frac{9}{5}|c|)^{2/3})$.
 The Jacobian matrix of $\bar{X}$ at $m_0$ is 
 $$\begin{pmatrix}
-\frac{1}{3}\bar{s}_0 & 0 & 0\\
* & \frac{5}{3}\bar{s}_0 & 1\\
* & \frac{2}{3}\bar{\zeta}_0 & \frac{2}{3}\bar{s}_0
 \end{pmatrix}$$
 giving the eigenvalue $-\frac{1}{3}\bar{s}_0$ in the direction of $R$. On $\{R=0\}$ 
 we get the two conjugate eigenvalues $\bar{s}_0(\frac{7}{6}\pm\frac{\sqrt{11}}{6}i).$ Thus $m_0$ is
  an hyperbolic equilibrium point, if $c>0$ (implying $\bar{s}_0>0$), it has one dimension stable manifold transverse to $S^2_+$,
   and a two dimensional unstable one, contained in $S^2_+$. If $c<0$, the situation is symmetric.
 
 \subparagraph{The chart $(ii)$.} Along $\p S^2_+$ we set 
 $\begin{cases}
\bar{\zeta}=\cos\omega\\
\bar{s}=\sin\omega
 \end{cases}$ 
 and proceed to the blow up $\rho=R^3\bar{\rho},\; s=R\sin\omega,\; \zeta=R^2\cos\omega$. 
 The pulled-back dynamics is 
 \begin{equation}
Y:\begin{cases}
R'=\frac{R}{1+\cos^2\omega}(\sin\omega(\cos\omega+\frac{1}{2}\sin^2\omega)+c\bar{\rho}\cos\omega)+O(R^2)\\
\omega'=\frac{1}{1+\cos^2\omega}(\cos\omega(2\cos\omega+\sin^2\omega)-c\bar{\rho}\sin\omega)+O(\bar{\rho}R)\\
\bar{\rho}'=-\frac{\bar{\rho}}{1+\cos^2\omega}(\sin\omega(1+\cos^2\omega+\cos\omega+1/2\sin^2\omega)+c\bar{\rho}\cos\omega)+O(\bar{\rho}R)
\end{cases}
 \end{equation}
 and is equivalent to 
  \begin{equation}
 \bar{Y}:\begin{cases}
 R'=R(\sin\omega(\cos\omega+\frac{1}{2}\sin^2\omega)+c\bar{\rho}\cos\omega)+O(R^2)\\
 \omega'=\cos\omega(2\cos\omega+\sin^2\omega)-c\bar{\rho}\sin\omega+O(\bar{\rho}R)\\
 \bar{\rho}'=-\bar{\rho}(\sin\omega(1+\cos^2\omega+\cos\omega+1/2\sin^2\omega)+c\bar{\rho}\cos\omega)+O(\bar{\rho}R)
 \end{cases}
 \end{equation}
In restriction to $\{\bar{\rho}=R=0\}$, we obtain 4 equilibrium points, namely, the solutions of
\linebreak $\cos\omega(\sin^2\omega+2\cos\omega)=0$. 
In addition to the trivial $\pm\pi/2$, we end up with
  $\cos\omega=1-\sqrt{2}$, this last equation gives two solutions 
  $\omega_0\in]\pi/2,\pi[$ and $-\omega_0$.
   All this zeros are simple (in the direction of $\omega$) so the d
   ynamics on $\p S^2_+$ can be
    deduced by the sign of $\cos\omega(2\cos\omega+\sin^2\omega)-
    c\bar{\rho}\sin\omega+O(\bar{\rho}R)$ on $\{\bar{\rho}=\omega=R=0\}$, 
    which is non negative.  
Actually, from $\omega_{0}$ and $-\omega_0$ we get two lines of zero in 
the plane $\{\bar{\rho}=0\}$,
  which are the blow up of the parabola $\zeta=-s^2/2$ (corresponding $\Gamma$).  

\begin{figure}[t]
	\centering
	\includegraphics[width=0.7\linewidth]{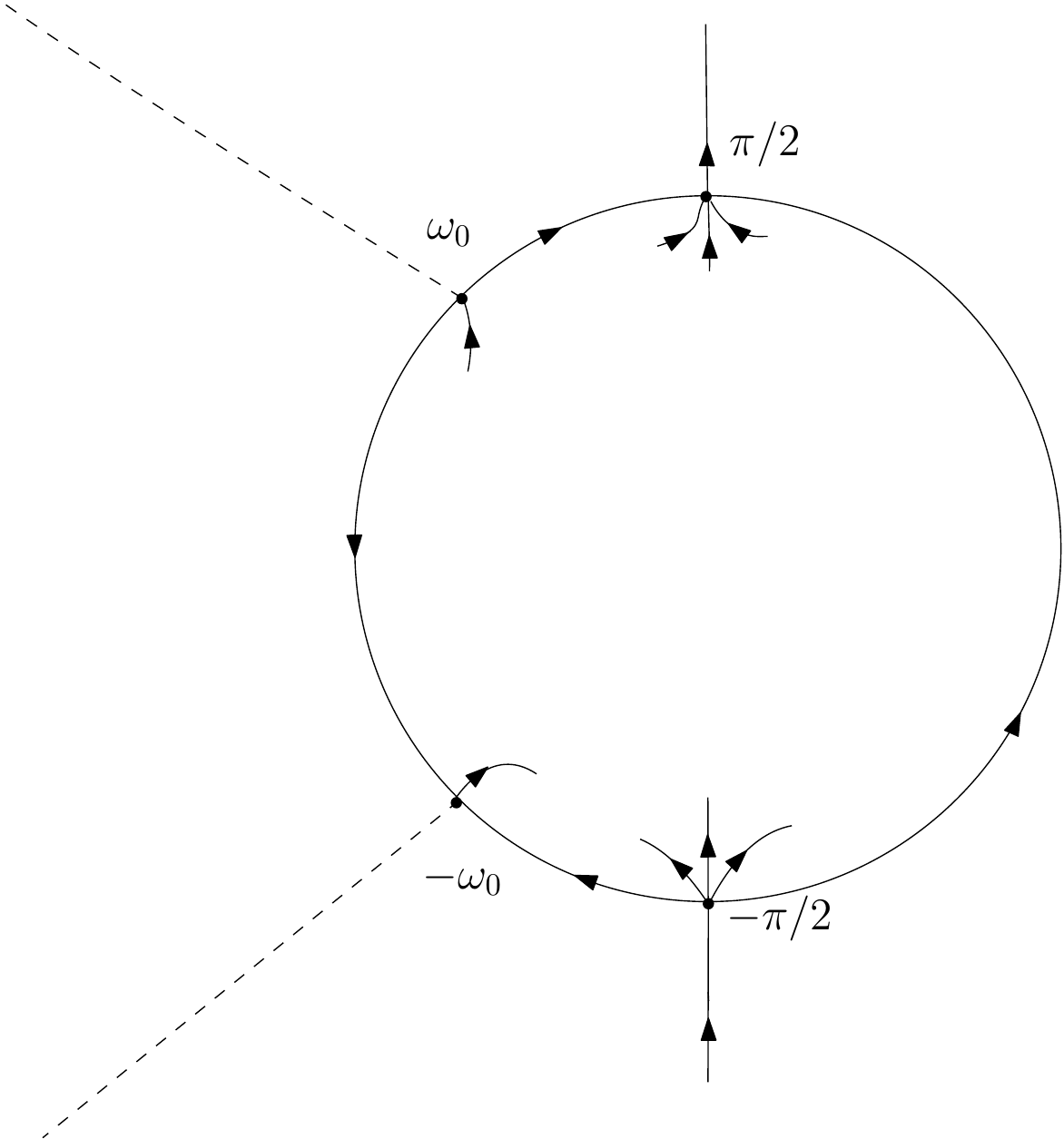}
	\caption{Phase portrait around $\p S^2_+$}
	\label{fig1}
\end{figure}

Let us write the Jacobian matrix of $\bar{Y}$ in the plane $\bar{\rho}=0$: 
$$J=\begin{pmatrix}
U(\omega) & * & *\\
0 & \Omega(\omega) & *\\
0 & 0 & \Gamma(\omega)
\end{pmatrix}$$ 
with $U(\omega)=\frac{\sin\omega}{2}(2\cos\omega+\sin^2 \omega)$, $\Omega(\omega)=
-\sin\omega(2\cos\omega+\sin^2 \omega-2\cos^2 \omega+2\cos\omega$, and
 $\Gamma(\omega)=-\sin\omega(\cos\omega+\frac{3}{2}\sin^2 \omega+2\cos^2 \omega)$. 
 We still need two informations:\newline 
- the eigenvalues of $\pm\pi/2$ in the radial direction, given by $U(\pm\pi/2)=\mp1$. \newline
- the eigenvalues of the 4 equilibria in the direction of $\bar{\rho}$, given by 
$R(\pm\pi/2)=\mp\frac{3}{2}$ and $R(\omega_0)<0$, $R(-\omega_0)>0$.
 Now we have a clear description of the phase portrait in a neighborhood of
  $\p S^2_+$ in Figure \ref{fig1}. $\pm\pi/2$ are hyperbolic saddles and
   $\pm\omega_0$ are hyperbolic in restriction to $S^2_+$ (but not in dimension $3$).
  We will see that the dynamics of $(\ref{hot})$ is also stable by perturbation by higher order terms. 

\paragraph{Global dynamics.} We restrict ourselves to the case $c>0$, the case $c<0$ being symmetric. We are now going to glue the studies in both charts to
 obtain the phase portrait on a whole neighborhood of the hemisphere.
 The main tool in that regard will be the following theorem from Poincaré and Bendixson, \cite{PB}. 

\begin{theorem}[Poincaré-Bendixson]
Let $X$ be a vector field in the plane, any maximal solution of $\dot{x}=X(x)$
 contained in a compact set, is either converging to an equilibrium point, a limit cycle, or a graphic, i.e., a close invariant curve union of a finite number of orbits connecting equilibria.
\end{theorem}

The equilibria of the flow restricted to $S^2_+$ (i.e., $R=0$) are 
as followed: $\pi/2\in S^1\cong\p S^2_+$ is a stable node, likewise, $-\pi/2$
 is an unstable node. The equilibrium $m_0$ 
is an unstable focus. $\omega_0$ and $-\omega_0$ are saddles:
 the stable manifold of $\omega_0$ (separatrix) and the unstable manifold of $-\omega_0$ are transverse to $S^2_+$ which contains their other invariant manifolds. Its unidimensional 
 unstable manifold is on $\partial S^2_+$.
Besides, for $\omega_0$, the opposite happens: It has a one dimensional 
stable manifold, and its unstable manifold is 
along $\p S^2_+$. 
Now we will prove that $\bar{X}$ does not have any periodic orbits. 
This, according to Poincaré-Bendixson, will allow us to link 
the trajectory coming from unstable directions to the stable manifolds
 belonging to other singular points in $S^2_+$.

\begin{figure}[t]
	\centering
	\includegraphics[width=0.7\linewidth]{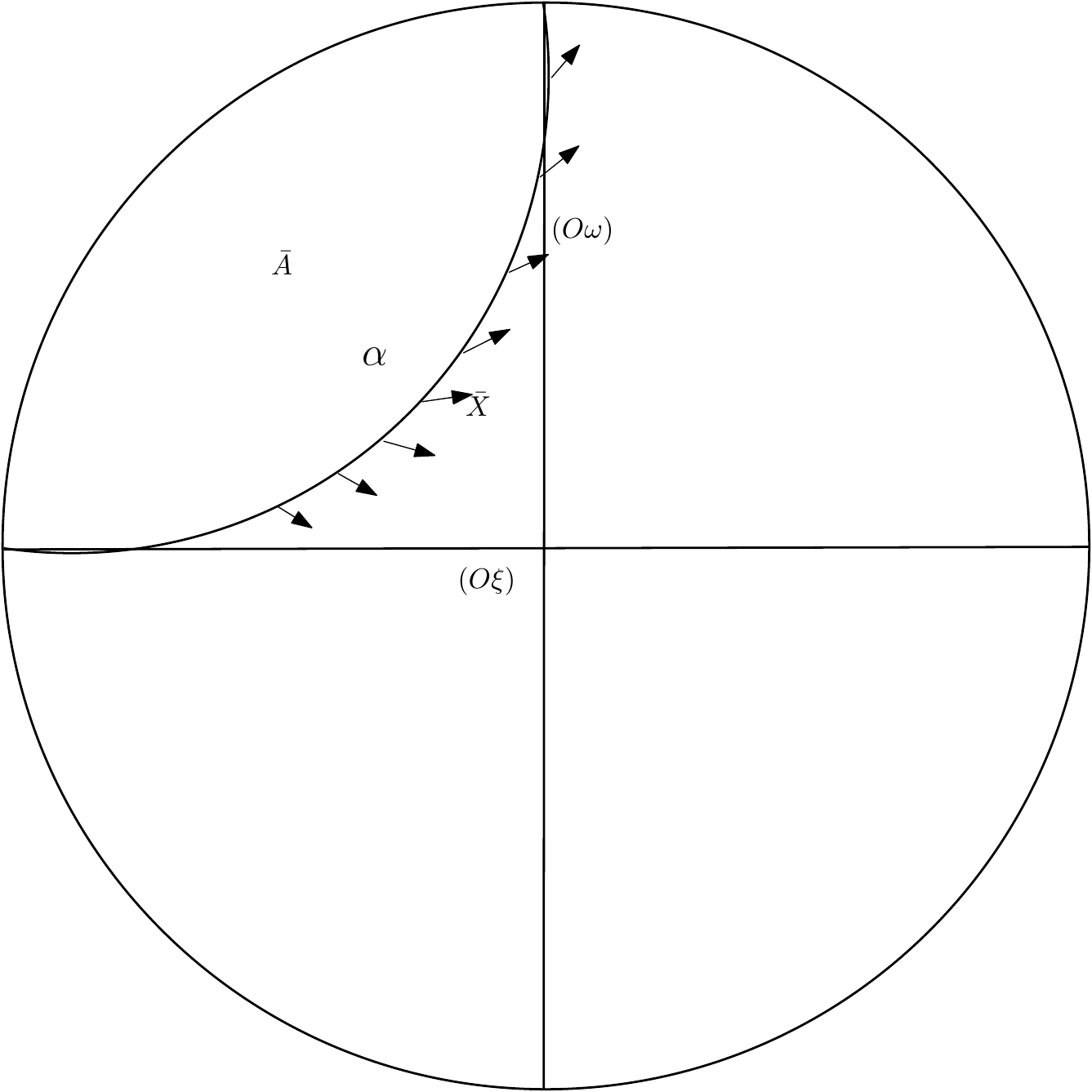}
\caption{Building $\bar{A}$}
\label{fig2}
\end{figure}

\begin{lemma}
$\bar{X}$ does not have a  periodic orbit on $S^2_+$. 
\end{lemma}
\noindent\textbf{Proof.} Between the charts $(i)$ and $(ii)$, we have the 
following change of coordinates : 
\begin{equation*}
\begin{cases}
\bar{s}=\frac{\sin\omega}{\bar{\rho}^{1/3}}\\
\bar{\zeta}=\frac{\cos\omega}{\bar{\rho}^{2/3}}.
\end{cases}
\end{equation*}
We can now define the two orthogonal axis $(O\bar{\zeta})$ and $(O\bar{s})$
 in $S^2_+$, $\p S^2_+$ included. In the chart $(ii)$, $(O\bar{\zeta})$ is going 
 from $\omega=\pi$ to $\omega=0$. Consider the convex domain $A$, such 
 that $\p A=(O\bar{s})_+\cup]\pi/2,\pi[\cup(O\bar{\zeta})_-$, then $m_0\in\text{Int}A$
  is the only equilibrium of $\bar{X}$ in $S=\text{Int}S^2_+$. The field $\bar{X}$ is 
  positively collinear to $(0\bar{\zeta})_+$ in $(0,0)$, and transverse to those axis 
  everywhere else. Thus, we can smooth the boundary of $A$ corresponding to 
  the part $(O\bar{\zeta})_-\cup(O\bar{s})_+$ by a curve $\alpha$ in order to make
   $\bar{X}$ transverse to $\p A$. See figure \ref{fig2}.

Denote $\bar{A}$ the part of $S^2_+$ such that $\bar{A}\subset A$ and
 $\p \bar{A}=]\pi/2;\pi[\cup\alpha$. Now $X$ is transverse to $\bar{A}$ and pointing
  outside $\bar{A}$. In $\bar{A}$, we have div$(X)=2\bar{s}>0$.
Now assume $\gamma$ is a periodic orbit of $\bar{X}$. By Jordan's theorem, 
$\gamma$ is the boundary of a compact set $D\subset S^2_+$, diffeomorphic to a disk. 
The result is then a consequence of the Poincaré-Hopf formula:
\begin{theorem}[Poincaré-Hopf]
Let $M$ be a compact manifold, and $X$ a vector field that has isolated zeros on $M$. 
Then $\sum_{i=1}^{m}\text{Index}(x_i)=\chi(M)$, where the $x_i$ are all the zeros of 
$X$ in $M$, and $\chi$ denotes the Euler characteristic.
\end{theorem}
$D$ being contractile, $\chi(D)=1$, hence $D$ contains at least one equilibrium 
point, and since $m_0$ is the only one in $S$, $m_0\in D$. 
As a result, either $\gamma$ lies in $\bar{A}$ or intersects $\alpha$. 
Let us consider the first alternative: $\gamma\subset\bar{A}$. We have 
\[0<\int_{D}\text{div}(\bar{X})\d\bar{\zeta}\wedge \d\bar{s}=
\int_{D} d(\iota_{\bar{X}}(\d\bar{\zeta}\wedge \d\bar{s}))=
\int_{\gamma}\iota_{\bar{X}}(\d\bar{\zeta}\wedge \d\bar{s})=0\]
by Stokes formula, which excludes that case. Now, note that all intersection
 points between $\gamma$ and $\alpha$ are transverse, since $\bar{X}$ is 
 transverse to $\alpha$: Thus, there is no tangency, and $\gamma$ intersects
  $\alpha$ twice. But this is also
excluded because $\bar{X}$ being transverse, it is only pointing outside $\bar{A}$. 
 $\Fin{\Box}$
\\

\begin{figure}
	\centering
	\includegraphics[width=1\linewidth]{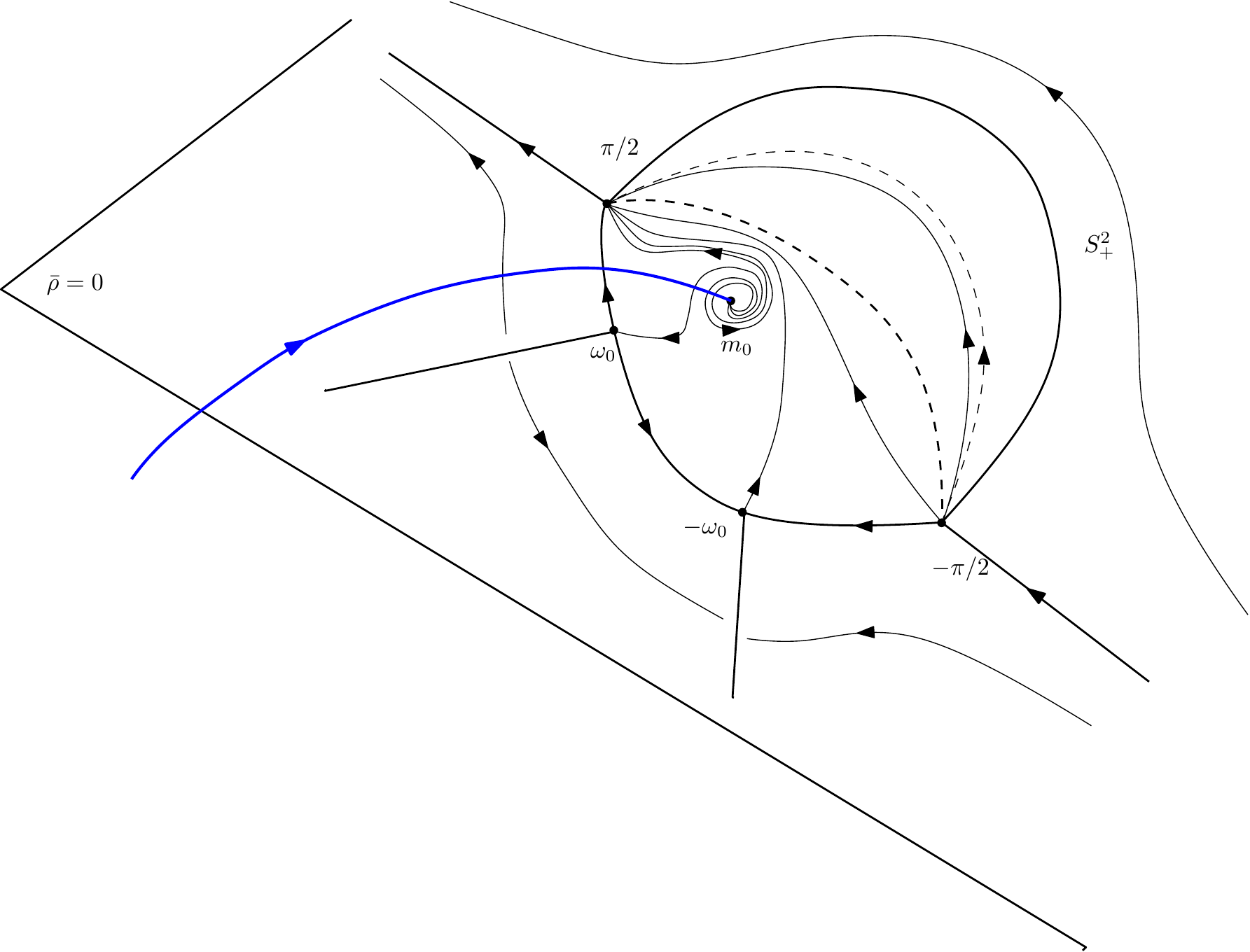}
	\caption{Phase portrait around $S^2_+$}
\label{fig3}
\end{figure}

By the Poincaré-Bendixson theorem above,
since there is no periodic orbits, in $\text{Int}\bar{A}$ every trajectory converges to
$m_0$ when the time tends to $-\infty$. $\omega_0\in\p\bar{A}$ has a stable manifold 
of dimension one, and this stable manifold lies inside $\text{Int}\bar{A}$ (at least close to $\omega_0$).
This implies that the stable manifold from $\omega_0$ converges to $m_0$ in negative infinite time. 
Apart from the equilibrium $\pi/2$, it is the only stable direction in $S$. That means all the other 
trajectories converge to the stable node (restricted to $S^2_+$) $\pi/2$, leading to the phase portrait of figure \ref{fig3}. As shown in \cite{orieux:2018:singularities}, the stable stratum of theorem 1, $S^s$, is the disjoint union of one dimensional stable manifolds to the equilibria lines from $\omega_0$. This submanifold expends on the critical locus and join $m_0$. The blow down leads to the local-global picture of figure \ref{fig0}

\paragraph{Back to the original system} 
The initial problem lies in dimension $k+2$ ($k=6$ in our initial control affine problem). From (\ref{0}), we see that when $\rho=0$, the $\xi$-component of (\ref{0}) vanishes. Thus, in the blown up coordinates, when $R=0$ (on $S^2_+$) or when $\bar{\rho}=0$, the spaces $\{\tilde{\xi}_2=\text{const}, \dots, \tilde{\xi}_k=\text{const}\}$ are preserved. Let us write their dynamics in the chart $(i)$ (with obvious notation with respect to ($\ref{0}$)): 
\begin{equation*}
\begin{cases}
\tilde{\xi}_2'=\bar{\rho}R^2\tilde{h}_2(R,\bar{s},\tilde{\xi})\\
\vdots\\
\tilde{\xi}_k'=\bar{\rho}R^2\tilde{h}_k(R,\bar{s},\tilde{\xi}).
\end{cases}
\end{equation*}
The blown up space is $S^2_+\times\R^{k-1}$, and the linear part of the total dynamics is the same as in (\ref{charti}), completed with zeros to obtain a $k+2$ matrix. 
As a result, in the initial system, the hyperbolic equilibrium point $m_0$ is replaced by a $k-1$ manifold of equilibria, denoted $N_0$, parametrized by $(\tilde{\xi}_2,\dots,\tilde{\xi_k})$. Each of these points have a stable manifold of dimension one in the direction of $R$, when $c>0$, (resp. unstable hen $c<0$): there exists a trajectory of (\ref{0}) converging to each of these points. 
The trajectories from the stable manifolds of $m_0=m_0(\tilde{\xi}_2,\dots,\tilde{\xi_k})$ are the one in theorem \ref{switchnilp}.

It remains to show 
that the trajectory coming from the stable manifold to $m_0$ is actually 
going to $m_0$ in finite time, for the original time $t$. We have been doing the following 
changes of times: $dt=\rho dt_1$, $dt_2=rdt_1$, $dt_3=Rdt_2$ ($t_3$ 
is the time in which we study the blown up system (\ref{charti})), so that $dt=\frac{\rho}{rR}dt_3$.
We will show that the interval of time 
from a point of the stable manifold to $m_0$ is finite. 
Assumption (\ref{eq:A}) implies among other things: $\rho=0\Rightarrow r>0$, 
so that in a neighborhood $O$ of $\Sigma$, we have $r>0$. Then in $O$, 
$r$ is bounded below and above by two positive constant $A>r>B>0$. 
In the blown up coordinates, $\rho=R^3\bar{\rho}$, so that the 
previously mentioned interval of time is 
$$\Delta t=\int_{t_3^0}^{+\infty}\frac{R^2(t_3)\bar{\rho}}{r(t_3)}dt_3<
\frac{1}{B}\int_{t_3^0}^{+\infty}R^2(t_3)\bar{\rho}(t_3)dt_3.$$ 
Notice that $\bar{\rho}$ is bounded by above by a positive constant 
$K$ along the trajectory in the stable manifold, since it converges to $m_0$. 

Now, the first line of system (\ref{charti}) is $R'=-\frac{1}{3} R\bar{s}+R^2v(z)$, with $v$ a smooth function.  Since $m_0\in\{\bar{s}>0\}$, if
$O$ is small enough, $R'<-CR$ since $v$ is bounded in $O$, for a constant $C>0$. 
Then as along the stable manifold to $m_0$, we have $R(t_3)<R_0e^{-Ct_3}$ by integration 
between a time $t_3$ and $t_3^0$. So that finally, 
$$\Delta t <K/B\int_{t_3^0}^{+\infty}R_0^2e^{-2Ct_3}dt_3<+\infty. $$
So $\bar{z}$ is reached in finite time. From figure \ref{fig3}, and the fact that
the $\tilde{\xi}_i$'s, $i>1$ are constant on $S^2_+\times\mathbb{R}^{k-2}$ and when $\bar{\rho}$ vanishes,
one can make the same time estimates to prove that the extremal
goes out of $S^2_+$ in finite time. Those extremals can be connected to the singular flow.

\subsubsection{Proof of theorem \ref{flotnilpreg}}

In the process of proving theorem \ref{switchnilp}, we obtained a good description of the singular flow around a point of $\Sigma_0$. 
We will make the proof when $c>0$, the opposite case 
being similar. 
The manifold $N_0=\{m_0\}\times\R^{k-1}_{\tilde{\xi}}\cap O_{\bar{z}}$ in the
blown up space is normally hyperbolic since $m_0$ is hyperbolic and 
each point is an equilibrium. In the direction of $R$, there 
is at one dimensional stable manifold at each of those points.
Thus, as in \cite{orieux:2018:singularities}, we can define
the global $k$-dimensional (or of codimension 2) stable manifold 
$$S^s_0=\underset{z\in N_0}{\bigcup} W^s(z),$$
and 
from \cite{hirsch2006invariant}, it is $\mathcal{C}^\infty$-smooth.
The construction of the strata $S^s$ is similar as, and detailed in \cite{orieux:2018:singularities}. 
So, the neighborhood $O_{\bar{z}}$ is stratified as wanted. The flow is trivially smooth
on $S_0$, it has also been proven that it is smooth on $S^s$. The only remaining thing to
prove is the smoothness on $S^s_0$. Define the contact time with $\Sigma_0$,
$t_{\Sigma_0}(z_0)=\int_{0}^{\infty}  \rho(t_1,z_0)dt_1$ for $z_0\in S^s_0$. One can then define the map $z_0\mapsto z_{\Sigma_0}(z_0)\in N_0$, the contact point with $N_0$ of the extremal from $z_0$. Since $S^s_0$ is a smooth foliation with leaves $W^s(z)$, this map is smooth. The flow is smooth as well on $D_-:=\{(t,z_0)/\ z_0\in S^s_0,\; t<t_{\Sigma_0}(z_0)\}$ 
because the singular locus has not been crossed yet. 
Restricted to $S^s_0$, the vector field has a submanifolds of zeros on $N_0$ and is tangent to the leaves $W^s(z),\; z\in N_0$, each on these zeros have only one nonzero eigenvalue: there is no resonance on the $W^s(z)$ and the system restricted to them is conjugated to its linear part. Now, from the smoothness of the foliation, the diffeomorphism realizing the conjugation depends smoothly on the parameters $(\tilde{\xi}_2,\dots,\tilde{\xi}_k)$. I.e. on $S^s_0$ the dynamics is conjugated to
$$\begin{cases}
R'=-1/3\bar{s_0}R\\
\bar{s}'=0\\
\bar{\zeta}'=0=\tilde{\xi}_2'=\dots=\tilde{\xi_k}'
\end{cases}
$$
Thus restricted to $S^s_0$, $R(t_3)=R_0e^{-1/3\bar{s}_0t_3}$ and the time of contact with the singular locus $\Sigma_0$ can be express as
$t_{\Sigma_0}(z_0)=\int_{0}^{+\infty}\frac{R^2(t_3)\bar{\rho}(t_3)}{r(\tilde{\xi}(t_3))}dt_3=3R_0^2\bar{\rho}_0/(2r_0)$ which is a smooth function of $z_0$.
The singular flow is smooth, from any point of $\Sigma_0$ according to proposition \ref{sing}. As a result, the map
$(t,z)\in D_+=\{(t,z_0)/\ z_0\in S^s_0,\; t>t_{\Sigma_0}(z_0)\}\mapsto z(t-t_{\Sigma_0}(z_0),z_{\Sigma_0}(z_0))=z(t,z_0)$ is also smooth, and theorem \ref{flotnilpreg} is proved.
The continuity is obtained by the same proof than in the $\Sigma_-$ case,
see \cite{orieux:2018:singularities}. 

\subsubsection{Proof of proposition \ref*{controlcont}}
Depending on the sign of $H_{12}$, the control does not have the same regularity.
In the coordinates of system $(\ref{syst1})$, when $t<\bar{t}$, $u(t)=(\cos\theta(t),\sin\theta(t))$, but from proposition \ref{sing}, when $t>\bar{t}$, $u(t)=u_s(t)=\frac{(-H_{02},H_{01})}{H_{12}}=\frac{r(-\sin\phi,\cos\phi)}{H_{12}}=\frac{r}{H_{12}}(\cos(\phi+\pi/2),\sin(\phi+\pi/2))$.  In the first alternative, the extremal reaches the singular locus at 
the equilibrium point in the time $\bar{t}=\int_0^{\infty}\rho dt_1$, and
we have $\theta(\bar{t})-\phi(\bar{t})=\pi/2$:  the control is continuous when 
the connection with the singular flow occurs. 
In the second one, $\theta(\bar{t})-\phi(\bar{t})=-\pi/2$, so that $\theta(\bar{t}^-)=\theta(\bar{t}^+)+\pi$. $\Fin{\Box}$

\begin{remark}
	\emph{From the phase portrait of figure \ref{fig3}, we can actually retrieve all three cases. Indeed, one can make a change of coordinates to integrate the parameters $\alpha$. More precisely, set $\tilde{\xi}_1=a(\xi)-1$. Furthermore, the cases $\Sigma_-$ can be seen as the West part of the phase portrait above the sphere $S^2_+$, the two lines of zeros corresponding to the partially hyperbolic equilibrium of \cite{orieux:2018:singularities}, to retrieve the global phase portrait one has to quotient the $s$ axis to keep $s$ in $\mathbb{S}^1$. The Est part above $S_+^2$ being the $\Sigma_+$ case. The dynamics is actually structurally stable, and the whole situation is contained in the nilpotent case $\Sigma_0$. }
\end{remark}

The following example is a twist of the nilpotent approximation
of the minimum time Kepler problem proposed in \cite{caillau2012minimum}.

\textit{Example.}
Let us exhibit a control-affine system with the kind of trajectory describe in theorem
\ref{switchnilp} when the final time is minimized.

Consider 
	\begin{equation}
\begin{cases}
\dot{x}(t)=F_0(x(t))+u_1(t)F_1(x(t))+u_2(t)F_2(x(t)), \quad
t\in[0, t_f], \quad u \in U\\
x(0)=x_0 \\
x(t_f)=x_f\\
t_f\rightarrow \min.
\end{cases}
\end{equation}
on $\R^4$
with 
\[\begin{cases}
F_0(x)=x_1\frac{\p}{\p x_3}+x_2\frac{\p}{\p x_4},\\ 
F_1(x)=x_2\frac{\p}{\p x_1}+\frac{\p}{\p x_3},\\
F_2(x)=\frac{\p}{\p x_2}.
\end{cases}\]
Then 
$$\text{rank}(F_{1}(x),F_2(x),F_{01}(x),F_{02}(x))=4,\; \forall x\in\R^4\setminus\{x_2=0\}.$$
The maximized Hamiltonian is $H^{\max}(x,p)=p_3x_1+p_4x_2+\sqrt{(p_1x_2+p_3)^2+p_2^2}$ and we have
\begin{equation}
\begin{cases}
\dot{x}_1=\frac{(p_1x_2+p_3)x_2}{\sqrt{((p_1x_2+p_3)^2+p_2^2)}}, & \dot{x}_3=x_1,\\
\dot{x}_2=\frac{p_2}{\sqrt{((p_1x_2+p_3)^2+p_2^2)}}, & \dot{x}_4=x_2.
\end{cases}
\end{equation}
The coordinates $x_3$ and $x_4$ are cyclic, so $p_3$ and $p_4$ are 
constant. Denote $p_3=-a$, $p_4=-c$, we get $p_1(t)=at+b$,
$p_2(t)=ct+d$, where $b=p_1(0),$ $d=p_2(0)$. Eventually: 
\begin{equation}
\dot{x}_2=\frac{ct+d}{\sqrt{((at+b)x_2-a)^2+(ct+d)^2}}.
\label{ex_2} 
\end{equation}
We also have $\Sigma=\{p_2=p_1x_2+p_3=0\}$, and the condition 
$H_{01}^2+H_{02}^2=H_{12}^2$ gives $\Sigma_0=\Sigma\cap\{p_1^2=p_3^2x_2^2+p_4^2\}.$
The contact time with $\Sigma$ has to be $\bar{t}=-\frac{d}{c}$. 
At $\bar{t}$, we must have $x_2(\bar{t})=-p_3/p_1(\bar{t})=-\frac{ac}{ad-bc}$. 
In order to reach $\Sigma_0$, we shall have 
$x_2^2(\bar{t})=\frac{1}{a^2}(p_1^2(\bar{t})-p_4^2)=\frac{(ad-bc)^2-c^2}{a^2c^2}:=\bar{x}$.
This gives an equation on the initial conditions $(a,b,c,d)$: 
\begin{equation}
(ad-bc)^2[(ad-bc)^2-c^2]-a^4c^4=0,
\label{cond0}
\end{equation}
choosing $a$ and $c$ non-zero. This condition imposes $z(\bar{t},z_0)\in\Sigma\Rightarrow z(\bar{t},z_0)\in\Sigma_0$. 
Now note that $x_2$ verifies a real ordinary differential equation (though, time dependent) $\dot{x}_2=f(t,x_2)$ with $f$ defined by (\ref{ex_2}).
$f$ is regular on $\R^2\setminus\{(\bar{t},\bar{x})\}.$ To regularize the dynamics of 
$x_2$ set $\d t=\sqrt{((at+b)x_2-a)^2+(ct+d)^2}\d s$ to obtain a continuous
dynamical system
in the plane:
\begin{equation*}
\begin{cases}
x_2'=ct+d\\
t'=\sqrt{((at+b)x_2-a)^2+(ct+d)^2}
\end{cases}
\end{equation*}
$(\bar{x},\bar{t})$ is its only equilibrium. Outside of it, $t'>0$. Choosing $c>0$, 
there exists a one dimensional stable manifold going to $(\bar{x},\bar{t})$, and thus a
trajectory converging to it in infinite time $s$. This implies the existence of a trajectory
for (\ref{ex_2}) such that $x_2(\bar{t})=\bar{x}$. Hence, together with condition 
(\ref{cond0}), there exists an extremal reaching $\Sigma_0$.



\begin{thebibliography}{99}
	

	
    \bibitem{agrachev04}
    Agrachev, A.~A.; Sachkov, Y.~L.
    {\em Control theory from the geometric viewpoint.} Springer, 2004.

	\bibitem{biolo}
	Agrachev, A. A.; Biolo, C.
	Switching in time-optimal problem: The 3D case with 2D control.
	{\em J.\ Dyn.\ Control Syst.}, \textbf{23} 2017, no.~3, 577--595.
	

	\bibitem{orieux:2018:singularities}
	J.-B. Caillau, J. F{\'e}joz, M. Orieux and R. Roussarie. 
	On singularities of min time control affine systems, (2018),
	\textit{submitted, 	arXiv:1907.02931} 
	
	
		
	\bibitem{dario}
	D. Barilari, Y. Chitour, F. Jean, D. Prandi, M. Sigalotti.
	On the regularity of abnormal minimizers for rank 2 sub-Riemannian structures, (2018),
	{\em  	arXiv:1804.00971} 
		
	\bibitem{caillau2012minimum}
	Caillau, J.-B.; Daoud, B.
	Minimum time control of the restricted three-body problem.
	{\em SIAM J.\ Control Optim.}\ \textbf{50} (2012), no.~6, 3178--3202.
	
\bibitem{intro}
Bonard, B.; Caillau, J.-B.
Introduction to non-linear optimal control.
{\em Advanced topics in control systems theory. Lecture Notes in Control and Inform. } Springer (2006).


	\bibitem{Durmotier}
	Dumortier, F. 
	Techniques in the Theory of Local Bifurcations: Blow-Up, Normal Forms, Nilpotent Bifurcations, Singular Perturbations, In: Schlomiuk D. (eds) Bifurcations and Periodic Orbits of Vector Fields
	{\em ATO ASI Series, vol 408. Springer, Dordrecht} (1993).

    \bibitem{PB}
    I. Bendixson. 
    Sur les courbes définies par des équations différentielles, {\em Acta Math., vol. 24, no 1, pp 1-88} (1901).
        
    \bibitem{hirsch2006invariant}
    Hirsch, M.~W.; Pugh, C.~C.; Shub, M.
    Invariant manifolds, Lecture Notes in Mathematics, vol. 583, Springer, Berlin Heidelberg, New York (1977).
    
    \bibitem{controlbox}
    Poggiolini, L. and Spadini, M.
    Bang-bang trajectories with a double switching time in the minimum time problem,
    {\em ENSAIM: COCV, vol. 22, 688 - 709} (2016).
	
\end{thebibliography}
\end{document}